\documentclass[12pt,a4paper]{article}

\usepackage{amsmath,amstext,amssymb,amscd}

\usepackage{graphicx}

\usepackage[russian,english]{babel}
\usepackage[cp1251]{inputenc}

\newcommand{\Xcomment}[1]{}

\newtheorem{theorem}{Theorem}[section]
\newtheorem{lemma}[theorem]{Lemma}

\newtheorem{prop}[theorem]{Proposition}

\newenvironment{proof}{\noindent{\bf Proof}\/}%
{\hfill$\qed$\medskip}

\def\qed{\Box}

\makeatletter \@addtoreset{equation}{section} \makeatother

\newenvironment{numitem1}{\refstepcounter{equation}\begin{enumerate}%
\item[(\thesection.\arabic{equation})]}{\end{enumerate}}

\newcommand{\refeq}[1]{(\ref{eq:#1})}  

 \makeatletter
\renewcommand{\section}{\@startsection{section}{1}{0pt}%
{-3.5ex plus -1ex minus -.2ex}{2.3ex plus .2ex}%
{\normalfont\Large}}
 \makeatother

 \makeatletter
\renewcommand{\subsection}{\@startsection{subsection}{2}{0pt}%
{-1.5ex plus -1ex minus -.2ex}{-1.5ex plus .2ex}%
{\normalfont\normalsize\bf}}
 \makeatother

 \newcommand{\SEC}[1]{\ref{sec:#1}}  
\newcommand{\SSEC}[1]{\ref{ssec:#1}}  

\def\Rset{{\mathbb R}}

\def\Zset{{\mathbb Z}}

\def\Ascr{{\cal A}}
\def\Bscr{{\cal B}}
\def\Cscr{{\cal C}}
\def\Dscr{{\cal D}}
\def\Escr{{\cal E}}

\def\Iscr{{\cal I}}

\def\Mscr{{\cal M}}
\def\Pscr{{\cal P}}

\def\tilde{\widetilde}
\def\hat{\widehat}
\def\bar{\overline}

\def\CXYXY{\Cscr_{X,Y,X',Y'}}

\begin{document}
\parskip=2pt

 \title{On universal quadratic inequalities for minors of TNN matrices}

 \author{Vladimir I.~Danilov\thanks{Central Institute of Economics and
Mathematics of the RAS, 47, Nakhimovskii Prospect, 117418 Moscow, Russia;
email: danilov@cemi.rssi.ru.}
 \and
Alexander V.~Karzanov
\thanks{Central Institute of Economics and Mathematics of
the RAS, 47, Nakhimovskii Prospect, 117418 Moscow, Russia; email:
akarzanov7@gmail.com.
}
  \and
Gleb A.~Koshevoy
\thanks{Institute for Information Transmission Problems of
the RAS, 19, Bol'shoi Karetnyi per., 127051 Moscow, Russia; email:
koshevoyga@gmail.com. }
  }
  
\date{}

 \maketitle

\vspace{-0.7cm}
 \begin{abstract}
For positive integers $n,n'$, we give a combinatorial characterization for the set of quadratic inequalities on minors that are valid for all $n\times n'$ totally nonnegative matrices. This is obtained as a consequence from our earlier results on stable quadratic identities on minors of matrices generated by flows in planar graphs via Lindstr\"om's construction.
 \medskip

{\em Keywords}\,: totally nonnegative matrix, minor, quadratic relation, Dodgson relation, planar graph, network flow

\medskip
{\em AMS Subject Classification}\, 05C75, 05E99

 \end{abstract}


\section{Introduction} \label{sec:intr}
 
In this paper we deal with real $n\times n'$ matrices in which all minors (viz. the determinants of square submatrices) are nonnegative, called \emph{totally nonnegative (TNN) matrices}. Our interest is focused on inequalities where each (left and right) side is formed by a sum of products of pairs of minors taken with unit coefficients. We refer to these as \emph{quadratic inequalities on minors} (\emph{QIM}) and say that such an inequality is \emph{universal} regarding the TNN matrices, or \emph{TNN-universal}, if it is valid for all $n\times n'$ TNN matrices. In literature many particular cases of such inequalities and equalities have been revealed, see e.g.~\cite{GK,Flm,FK,FV,muhl,Po,skan,SV} (where the list can be continued). Our goal is to give a complete characterization for the class of TNN-universal QIM. We state this in a combinatorial way and accompany by an ``efficient'' algorithm which, given an input inequality (encoded via the row and column index sets of involved minors), recognizes whether it is TNN-universal or not.

To obtain the desired characterization, we rely on constructions and results from our earlier work~\cite{DKK2}. The subject of that work differs from our present one in two points. First, it dealt with a wider scope of matrices, namely, the ones generated, via an analog of Lindstr\"om's construction~\cite{Li}, by systems of disjoint paths in a planar directed graph whose vertices take weights in an arbitrary commutative semiring. Following terminology from~\cite{DKK2}, we refer to such path systems as \emph{flows}, and to the obtained matrices as \emph{flow-generated} ones, or \emph{FG-matrices}. This is a representative class of matrices, which includes the TNN ones as a special case (cf.~\cite{Br}). Second, \cite{DKK2} studied quadratic \emph{identities}, rather than inequalities, on minors of FG-matrices. Popular examples in that area are formed by algebraic and tropical Pl\"ucker identities.

The main result in~\cite{DKK2} characterized the set of \emph{stable} quadratic identities on minors of FG-matrices, where the stability means that the identity is valid regardless of the planar graph and semiring generating an FG-matrix.
This is based on associating to each pair of minors involved in an examined identity $\Iscr$ a certain set of \emph{feasible planar matchings}. It was shown that $\Iscr$ is stable if and only if the families of such matchings in the left and right sides of $\Iscr$ are \emph{balanced}; namely, for each matching, the amounts of its occurrences in these sides are equal. (Note that subsequently the matching approach was extended in~\cite{DK} to obtain a characterization for the set of quadratic identities on minors of quantum matrices, which turned out to be flow-generated as well.)

This paper is organized as follows. Section~\SEC{theorem} contains basic definitions, including the formulation of quadratic inequalities on minors (QIM) of our interest. Then it recalls the notion of feasible (planar) matchings associated with an appropriate pair of (row-column index sets of) minors, and states our main result (Theorem~\ref{tm:univ_relat}) asserting that a QIM is TNN-universal if and only if certain inequalities on amounts of occurrences of feasible matchings in the left and right sides of the QIM are valid. Section~\SEC{back} describes a machinery of  constructions and tools from~\cite{DKK2} needed to us. Using these, we outline a proof of Theorem~\ref{tm:univ_relat} in Sect.~\SEC{proof}. The concluding Sect.~\SEC{illustr} demonstrates two appealing illustrations to our matching method, one of which is of submodular character, and the other generalizes Dodgson's type relations to TNN matrices. (It should be emphasized that in this paper we are focused mostly on a general problem of characterizing QIM for TNN matrices, and a few of particular cases are demonstrated only as illustrations.)


\section{Definitions and the theorem} \label{sec:theorem}

For positive integers $n,n'$, let $\Escr^{n,n'}$ denote the set of pairs $(I,J)$ such 
that $I\subseteq[n]$, \; $J\subseteq[n']$ and $|I|=|J|$, where $[n'']$ stands for $\{1,\ldots,n''\}$. For a real $n\times n'$ matrix $Q$ and $(I,J)\in\Escr^{n,n'}$, let $\Delta(I|J)=\Delta_Q(I|J)$ denote the determinant of the square submatrix of $Q$ with the row set $I$ and the column set $J$, called the $(I|J)$-\emph{minor} of $Q$. For some reasons, when $I,J=\emptyset$, we define $\Delta_Q(I|J)=1$. A matrix $Q$ is said to be \emph{totally nonnegative}, or a \emph{TNN matrix} for short, if $\Delta(I|J)\ge 0$ for all $(I,J)\in\Escr^{n,n'}$.

To introduce quadratic inequalities of our interest, consider disjoint subsets $X,Y$ of $[n]$ and disjoint subsets $X',Y'$ of $[n']$. Let $\CXYXY$ denote the set of pairs $(C\subseteq Y,\,C'\subseteq Y')$ such that
  \begin{equation} \label{eq:XCXC}
  |X\cup C|=|X'\cup C'| \quad\mbox{and}\quad |X\cup (Y-C)|=|X'\cup (Y'-C')|.
  \end{equation}
  
One can see that $\CXYXY\ne \emptyset$ if and only if $2|X|+|Y|=2|X'|+|Y'|$ (in particular, $|Y|+|Y'|$ is even). Also~\refeq{XCXC} implies $|Y|-|Y'|=2|C|-2|C'|$. For brevity we will write $\bar C$ for $Y-C$, and $\bar C'$ for $Y'-C'$ (assuming that $Y,Y'$ are fixed); also we will write $XC$ for $X\cup C$ (and so on); then the equalities in~\refeq{XCXC} can be rewritten as $|XC|=|X'C'|$ and $|X\bar C|=|X'\bar C'|$.

We call each pair $(C,C')\in \CXYXY$ \emph{proper} for $X,Y,X',Y'$ and consider two families $\Ascr$ and $\Bscr$ consisting of such pairs admitting multiplicities; for convenience we use notation $\Ascr\Subset \CXYXY$ and $\Bscr\Subset \CXYXY$. (Here $\Subset$ means that one can include in $\Ascr$ and $\Bscr$ several copies of the same element $(C,C')$ of $\CXYXY$. Multiplicities may appear when we wish to add up appropriate inequalities, as though forming a semigroup on them.) We are interested in pairs $(\Ascr,\Bscr)$ satisfying the inequality
  \begin{multline} \label{eq:q_relat}
  \sum\nolimits_{(A,A')\in \Ascr} \Delta_Q(XA|X'A')\, \Delta_Q(X\bar A|X'\bar A') \\
      \ge \sum\nolimits_{(B,B')\in \Bscr} \Delta_Q(XB|X'B')\, \Delta_Q(X\bar B|X'\bar B')\quad
   \end{multline}
for every $n\times n'$ TNN matrix $Q$; in this case we say that the pair $(\Ascr,\Bscr)$ generates a \emph{universal quadratic inequality} for TNN matrices, and that $(\Ascr,\Bscr)$ is \emph{TNN-universal}.

In this paper, relying on results on stable quadratic identities for minors of flow-generated matrices established in~\cite{DKK2}, we give a complete characterization of TNN-universal pairs $(\Ascr,\Bscr)$ (for any $n,n'$).


A characterization that we obtain involves a geometric-combinatorial construction described as follows. Assume that the elements of $Y$ (resp. $Y'$) are disposed on the lower (resp. upper) half of a circumference $O$ in the plane, so that the elements of each of $Y,Y'$ follow in the increasing order from left to right. For a proper pair $(C,C')$, let us refer to the elements of $C\cup C'$ (regarded as points in $O$) as \emph{white}, and to the elements of $\bar C\cup\bar C'$ as \emph{black}. 
  \medskip
  
\noindent\textbf{Definition.} Since $|Y|+|Y'|$ is even, one can partition $Y\sqcup Y'$ into 2-element subsets, or \emph{couples} (where $\sqcup$ denotes the disjoint union). We refer to a set $M$ of couples of this sort as a (perfect) \emph{matching} on $Y\sqcup Y'$ and say that $M$ is \emph{feasible} for $(C,C')$ if the following conditions hold:
  \begin{numitem1} \label{eq:feas_match}
  \begin{itemize}
\item[(i)] for each couple $\pi\in M$, if the elements of $\pi$ have different colors, then both elements lie in the same half of $O$ (either in the lower or in the upper one); whereas if they have the same color, then they lie in different halves;
\item[(ii)] $M$ is planar, in the sense that the chords to $O$ connecting the couples in $M$ are pairwise non-intersecting.
 \end{itemize}
 \end{numitem1}

\vspace{-0.2cm}
\begin{center}
\includegraphics[scale=0.6]{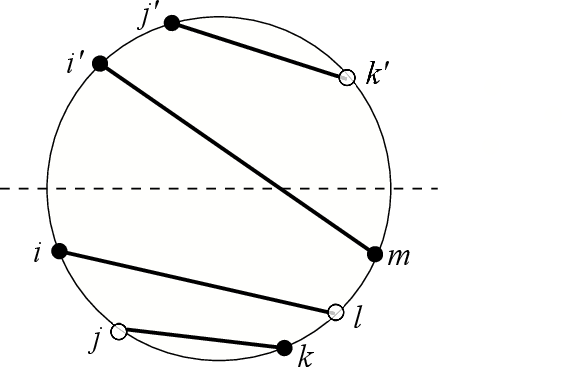}
\end{center}
\vspace{-0cm}
 
For convenience, to distinguish between elements of $Y$ and $Y'$ in illustrations, we denote the latter ones with primes. An example illustrating~\refeq{feas_match} is drawn in the above  picture. Here: $Y=(i<j<k<\ell<m)$,\, $Y'=(i'<j'<k')$,\, $C=\{j,\ell\}$,\, $C'=\{k'\}$,\, $\bar C=\{i,k,m\}$,\, $\bar C'=\{i',j'\}$, and a feasible matching is represented by 4 chords. 

A proper $(C,C')\in\CXYXY$ may admit one, two or more feasible matchings; we denote the set of these by $\Mscr_{C,C'}$. For $\Ascr\Subset \CXYXY$, let $\Mscr(\Ascr)$ denote the family of all feasible matchings (respecting multiplicities) occurring in the sets $\Mscr_{C,C'}$ over all $(C,C')\in\Ascr$. Also for a perfect matching $M$ on $Y\sqcup Y'$, denote the number of occurrences of this $M$ in $\Mscr(\Ascr)$ by $\#_M(\Ascr)$. Our criterion for TNN-universal QIMs is stated as follows.

\begin{theorem} \label{tm:univ_relat}
For $n,n',X,Y,X',Y'$ as above, let $\Ascr,\Bscr\Subset\CXYXY$. Then the following are equivalent:
 \begin{itemize}
\item[\rm(i)] the pair $(\Ascr,\Bscr)$ is TNN-universal, i.e.~\refeq{q_relat} is valid for all $n\times n'$ totally nonnegative matrices $Q$;
 \item[\rm(ii)] for each matching $M\in\Mscr(\Bscr)$, there holds $\#_M(\Ascr)\ge\#_M(\Bscr)$.
  \end{itemize}
  \end{theorem}
  

The proof of this theorem relies on results established in~\cite{DKK2} in context of stable quadratic identities for minors of flow-generated matrices, and on some technical tools elaborated there. In fact, \cite{DKK2} provides us with all ingredients sufficient to prove the theorem. We give a review in detail in Sect.~\SEC{back}, sometimes accompanying this with enlightening illustrations, and then outline a proof of Theorem~\ref{tm:univ_relat} in Sect.~\SEC{proof}. 

We finish this section with an example illustrating Theorem~\ref{tm:univ_relat} (this is a special case of a relation in~\cite[Theorem~B]{FV}).
\medskip

\noindent\textbf{Example}. Let $n=n'=5$, $X=X'=\emptyset$, $Y=[5]$, and $Y'=\{1',2',3',4',5'\}$. Define 
  $$
  A:=B:=12,\quad A'_1:=4'5',\quad A'_2:=3'4', \quad B':=3'5'.
  $$
Then $\bar A=\bar B=345$,\, $\bar A'_1=1'2'3'$,\, $\bar A'_2=1'2'5'$, and $\bar B'=1'2'4'$. (For brevity we write $ij$ for $\{i,j\}$, and $ijk$ for $\{i,j,k\}$.) Take the families $\Ascr:=\{(A,A'_1),(A,A'_2)\}$ and $\Bscr:=\{(B,B')\}$. Then relation~\refeq{q_relat} is specified as
  $$
  \Delta(12|4'5')\,\Delta(345|1'2'3')+\Delta(12|3'4')\,\Delta(345|1'2'5')
        \ge \Delta(12|3'5')\,\Delta(345|1'2'4')
  $$
(omitting matrix symbol $Q$). Relying on implication (ii)$\to$(i) in Theorem~\ref{tm:univ_relat}, in order to establish validity of the above inequality, it suffices to examine the corresponding feasible matchings. A routine verification shows that: 
  \begin{itemize}
\item $\Mscr_{A,A'_1}$ is formed by the unique matching $M_1:=\{14,23,3'4',2'5',51'\}$;
\item $\Mscr_{A,A'_2}$ is formed by two matchings $M_2:=\{14,23,2'3',4'5',51'\}$ and $M_3:=\{14,23,1'4',2'3',55'\}$.
 \item $\Mscr_{B,B'}$ is formed by two matchings, which are just $M_1$ and $M_2$.
    \end{itemize}
Then $\Mscr(\Ascr)=\{M_1,M_2,M_3\}$ and $\Mscr(\Bscr)=\{M_1,M_2\}$. This implies $\Mscr(\Ascr)\supset\Mscr(\Bscr)$, and the result follows.
\medskip

\noindent\textbf{Remark.} In light of a reduction to matchings figured in Theorem~\ref{tm:univ_relat}, it suffices to examine universal quadratic relations under the assumptions that $\min\{|X|,|X'|\}=0$, $Y=[m]$ and $Y'=[m']$, where $m=|Y|$ and $m':=|Y'|$. Then $\CXYXY$ consists of the pairs $(C\subseteq [m],\, C'\subseteq [m'])$ such that $|C|-|C'|=(m-m')/2$.
 \medskip


\section{Backgrounds}  \label{sec:back}

We start with a construction of flow-generated functions. As before, $n,n'\in\Zset_{>0}$.

\subsection{Planar networks and flows.} \label{ssec:flows}

By a \emph{path} in a directed graph we mean a sequence $P=(v_0,e_1,v_1,\ldots,e_k,v_k)$ where each $e_i$ is an edge connecting vertices $v_{i-1},v_i$. An edge $e_i$ is called \emph{forward} if it is directed from $v_{i-1}$ to $v_i$,
denoted as $e_i=(v_{i-1},v_i)$, and \emph{ backward} otherwise (when
$e_i=(v_i,v_{i-1})$). The path $P$ is called {\em directed} if it has no
backward edges, and {\em simple} if all vertices $v_i$ are different. 

By a \emph{planar network} we mean a finite directed graph $G=(V,E)$ such that: (a) $G$ contains two distinguished subsets of vertices $S=\{s_1,\ldots,s_n\}$ and $T=\{t_1,\ldots, t_{n'}\}$, called \emph{sources} and \emph{sinks}, respectively; (b) $G$ is \emph{acyclic} (has no directed cycles) and embedded in a compact convex region in the plane (admitting intersections of edges only in common vertices), and (c) the sources and sinks, also called \emph{terminals}, are disposed in the boundary $O$ of this region in the cyclic order $s_n,\ldots,s_1,t_1,\ldots,t_n$ clockwise. (We may think of $O$ as a ``circumference''.) It is useful to assume that $G$ has no \emph{redundant} edges, which means that any edge of $G$ belongs to a directed path from $S$ to $T$.

Let $(I,J)\in\Escr^{n,n'}$ (i.e. $I\subseteq[n]$, $J\subseteq[n']$ and $|I|=|J|$). By an $(I|J)$-\emph{flow} we mean a collection $\phi$ of \emph{pairwise} (vertex) \emph{disjoint} directed paths in $G$ going from the source set $S_I:=\{s_i\colon i\in I\}$ to the sink set $T_J:=\{t_j\colon j\in J\}$. (From the planarity of $G$ and the disposition of $S\cup T$ on $O$ it follows that the orders on $S_I$ and $T_J$ of the ends of paths in $\phi$ are agreeable, namely, if $I=(\alpha(1)<\cdots<\alpha(k))$ and $J=(\beta(1)<\cdots <\beta(k))$, then $i$-th path in $\phi$ goes from $s_{\alpha(i)}$ to $t_{\beta(i)}$.) 

We denote the collection of $(I|J)$-flows in $G$ by $\Phi_{I|J}=\Phi^G_{I|J}$.  This is related to $(I|J)$-minors of $n\times n'$ matrices due to Lindstr\"om's construction. More precisely, consider a nonnegative real-valued \emph{weighting} $w:V\to\Rset_+$ on the vertices. It generates the function $f_{G,w}$ on $\Escr^{n,n'}$ defined by
  \begin{equation} \label{eq:alg_f}
  f_{G,w}(I|J):=\sum_{\phi\in\Phi_{I|J}}\prod_{v\in V_\phi} w(v),
  \qquad (I,I')\in\Escr^{n,n'},
  \end{equation}
where $V_\phi$ is the set of vertices occurring in a flow $\phi$. (Note that $f_{G,w}(I|J)=0$ if $I,J\ne\emptyset$ and the set of $(I|J)$-flows in $G$ is empty. At the same time, when $I,J=\emptyset$, we prefer to put $f_{G,w}(I|J):=1$.) We refer to $f_{G,w}$ obtained in this way as a \emph{flow-generated function} (related to graph $G$ and weighting $w$). By Lindstr\"om's lemma~\cite{Li} (in the simplest version),
  \begin{numitem1} \label{eq:lind}
if $Q=Q_{G,w}$ is the $n\times n'$ matrix whose entries $q_{ij}$ are defined to be equal to  $f_{G,w}(i|j)=\sum_{\phi\in\Phi_{i|j}}\prod_{v\in V_\phi} w(v)$ (where for brevity we write $i|j$ for $\{i\}|\{j\}$), then for any $(I|J)\in\Escr^{n,n'}$, the minor value $\Delta_Q(I|J)$ is equal to $f_{G,w}(I|J)$.
 \end{numitem1}
(See also~\cite{GV}.) Since $w$ is nonnegative, all entries of $Q$, as well as all minor for $Q$, are nonnegative; so $Q$ is a TNN matrix. (As is said in~\cite{DKK2}, an analog of Lindstr\"om's lemma is valid for a wider class of weights $w$, namely, when $w$ takes values in an arbitrary commutative semiring; in particular, in a tropical one. For our purposes, it suffices to deal with nonnegative real weights only.) A converse property is valid as well. Namely, it is shown in Brente~\cite{Br} that 
  \begin{numitem1} \label{eq:brente}
any $n\times n'$ TNN matrix can be obtained by use of Lindstr\"om's construction, by choosing some planar network $G=(V,E)$ (subject to above conditions) and some nonnegative weighting $w$ on $V$.
 \end{numitem1}
 
Next, to simplify technical details, it is convenient to modify the network $G$ as follows. Split each vertex $v\in V$ into two vertices $v',v''$
(placing them in a small
neighborhood of $v$ in the plane) and connect them by edge $e_v=(v',v'')$,
called a \emph{split-edge}. Each edge $(u,v)$ of $G$ is replaced by an edge
going from $u''$ to $v'$. Also for each
$s_i\in S$, we add a new source $\hat s_i$ and edge $(\hat s_i,s'_i)$, and
for each $t_j\in T$, add a new sink $\hat t_j$ and edge $(t''_j,\hat t_j)$;
we refer to such edges as \emph{extra} ones. We denote the resulting graph by $\hat G=(\hat V,\hat E)$, and take $\hat S:=\{\hat s_1,\ldots,\hat s_n\}$ and $\hat T:=\{\hat t_1,\ldots,\hat t_{n'}\}$ as the sets of sources and sinks in it, respectively. The picture illustrates the transformation of $G$ into $\hat G$; here in the right fragment the split-edges are drawn bold, and the extra edges are horizontal and vertical.

\vspace{-0.2cm}
\begin{center}
\includegraphics[scale=0.8]{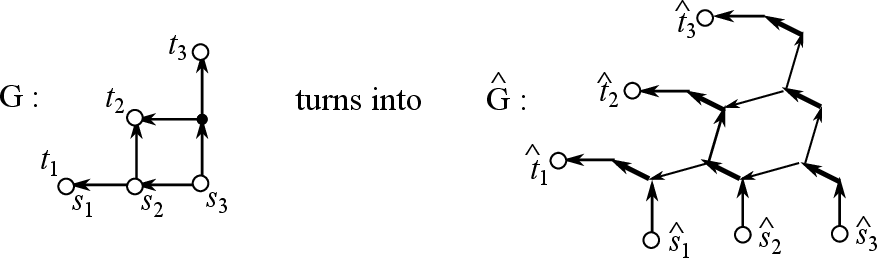}
\end{center}
\vspace{-0cm}

Clearly the modified graph $\hat G$ is acyclic, and the facts that the original graph $G$ is planar and has no redundant edges easily implies that (the induced embedding of) $\hat G$ is planar as well. For any $i\in[n]$ and $j\in[n']$, there is a natural one-to-one correspondence between the directed paths from $s_i$ to $t_j$ in $G$ and the ones from $\hat s_i$ to $\hat t_j$ in $\hat G$. This is extended to a one-to-one correspondence between flows, and for $(I,J)\in \Escr^{n,n'}$, we keep notation $\Phi_{I|J}$ for the set of flows in $\hat G$
going from $\hat S_I:=\{\hat s_i \,\colon i\in I\}$ to $\hat T_{J}:=\{\hat
t_j\,\colon j\in J\}$. It is convenient to transfer a weighting $w$ from the vertices $v$ of $G$ to the split-edges $e_v$ of ${\hat G}$, by setting $w(e_v):=w(v)$. Then the corresponding flows in both networks have equal weights (which are the products of the weights of vertices or split-edges in the flows). This implies that in both cases the functions on $\Escr^{n,n'}$ generated by corresponding flows are the same. Also $\hat G$ possesses the following useful properties:
  \begin{numitem1} \label{eq:prop_hatG}
\begin{itemize}
 \item[(i)] each non-terminal vertex is incident with exactly one split-edge, $e=(u,v)$ say; this $e$ is the only edge of $\tilde G$ leaving $u$ and the only edge entering $v$;
  \item[(ii)] each source  (sink) has exactly one leaving edge and no entering edge (resp. one entering edge and no leaving edge).
   \end{itemize}
   \end{numitem1}

\subsection{Double flows.} \label{ssec:double}

Let $(C,C')\in\CXYXY$. By a \emph{double flow} related to $(C,C')$ we mean a pair $(\phi,\phi')$ where $\phi$ is an $(XC|X'C')$-flow, and $\phi'$ an $(X\bar C|X'\bar C')$ flow in $\hat G$ (going from $\hat S_{XC}$ to $\hat T_{X'C'}$, and from  $\hat S_{X\bar C}$ to $\hat T_{X'\bar C'}$, respectively). We write $V_\phi$ and $E_\phi$ for the sets of vertices and edges occurring in $\phi$, respectively, and similarly for $\phi'$. We will rely on the next two lemmas from~\cite{DKK2} (giving short explanations below). Here $m:=|Y|$, $m':=|Y'|$, and we write $K\triangle L$ for the symmetric difference $(K-L)\cup (L-K)$ of subsets $K,L$ of a set.
\begin{lemma} \label{lm:Ephi-Ephip}
{\rm (cf. \cite[Lemma~2.1]{DKK2})}
For a double flow $(\phi,\phi')$, the set $E_\phi\triangle\, E_{\phi'}$ is partitioned into the edge sets of pairwise vertex disjoint circuits $C_1,\ldots,C_d$ (for some $d$) and simple paths $P_1,\ldots, P_p$ (where $p=(m+m')/2$), and each $P_i$ connects either $\hat S_{C}$ and $\hat S_{\bar C}$, or $\hat S_{C}$ and $\hat T_{C'}$, or $\hat S_{\bar C}$ and $\hat T_{\bar {C'}}$, or $\hat T_{C'}$ and $\hat T_{\bar {C'}}$. In each of these circuits and paths, the edges of $\phi$ and the edges of $\phi'$ have opposite directions (say, the former edges are forward, and the latter ones are backward).
\end{lemma}

\noindent (Here the facts that $E_\phi\triangle\, E_{\phi'}$ is partitioned into the edge sets of pairwise disjoint paths and cycles, and that neighboring edges in these paths/cycles contained in different $\phi$ and $\phi'$ have opposite directions easily follow from~\refeq{prop_hatG}(i). And~\refeq{prop_hatG}(ii) implies that any edge $e$ incident to $\hat S_{X}$ or $\hat T_{X'}$ belongs to both $\phi$ and $\phi'$; therefore, $e$ is not in $E_\phi\triangle\, E_{\phi'}$, implying that the paths $P_1,\ldots, P_p$ have their end vertices in $\hat S_{Y}\cup \hat T_{Y'}$, and that $p=(m+m')/2$. The remaining properties in the lemma follow from planarity.)
 \smallskip

\begin{figure}[htb]
\begin{center}
\vspace{-0.1cm}
\includegraphics[scale=0.9]{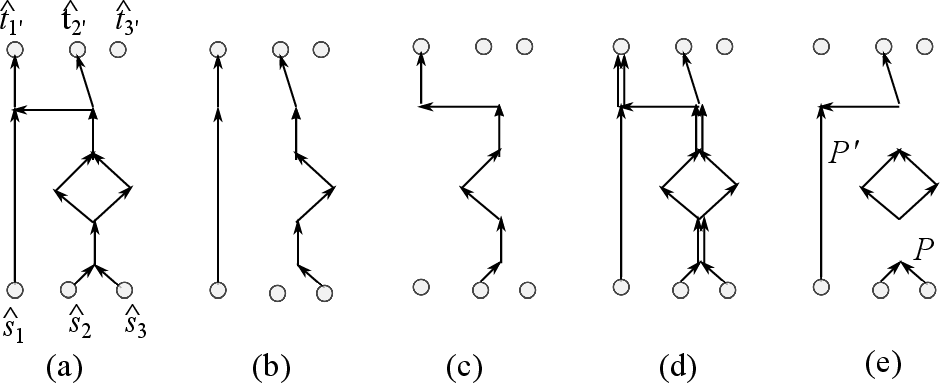}
\end{center}
\vspace{-0.5cm} 
\caption{(a) $\hat G$;\; (b) $\phi$;\; (c) $\phi'$;\; (d) $E_\phi\sqcup E_{\phi'}$;\; (e) $E_\phi\triangle E_{\phi'}$}
 \label{fig:double}
\end{figure}

This lemma is illustrated by the example shown in Fig.~\ref{fig:double}. Here $X=\emptyset$, $X'=\{1'\}$, $Y=\{1,2,3\}$, $Y'=\{2'\}$, $C=\{1,3\}$, $C'=\{2'\}$,  symbol $\sqcup$ means the union respecting multiplicities, and some edges and nonterminal vertices are hidden.

\subsection{Rearranging double flows.} \label{sec:rearrange}
Next we explain how to rearrange a double flow $(\phi,\phi')$ for $(C,C')$ so
as to obtain a double flow for another proper pair $(D,D')$ in $\CXYXY$. Let $\Pscr(\phi,\phi')$ denote the set of paths $P_1,\ldots,P_p$ as in Lemma~\ref{lm:Ephi-Ephip}, where $p=(m+m')/2$. For $P\in\Pscr(\phi,\phi')$, let $\pi(P)$ denote the pair of elements in $Y\sqcup Y'$ corresponding to the end vertices of $P$, and $E_P$ the set of edges of $P$. By Lemma~\ref{lm:Ephi-Ephip}, $\pi(P)$ belongs to either $C\times \bar C$ or $C\times C'$ or $C'\times \bar C'$ or $\bar C\times \bar C'$. Define
  $$
  M(\phi,\phi'):=\{\pi(P)\,\colon P\in\Pscr(\phi,\phi')\}.
 $$
This set forms a \emph{perfect matching} on $Y\sqcup Y'$. Moreover, since the paths in $\Pscr(\phi,\phi')$ are pairwise disjoint and connect vertices in the boundary of the corresponding convex region, $M(\phi,\phi')$ is planar. So $M(\phi,\phi')$ is a feasible matching for $(C,C')$ (cf.~\refeq{feas_match}).
   \begin{lemma} \label{lm:path_switch}
{\rm (cf. \cite[Lemma~2.2]{DKK2})}
Consider an arbitrary subset $M_0\subseteq M(\phi,\phi')$. Let 
  $$
  V_{M_0}:=\cup(\pi\in M_0)\;\; \mbox{and}\;\; E^0:=\cup(E_{P}\,\colon P\in\Pscr(\phi,\phi'),\,\pi(P)\in M_0).
  $$
Define $D:=C\triangle (V_{M_0}\cap Y)$ and $D':=C'\triangle (V_{M_0}\cap Y')$. Then there are a unique $(XD|X'D')$-flow $\psi$ and a unique $(X\bar D|X'\bar D')$-flow $\psi'$ such
that $E_\psi=E_\phi\triangle E^0$ and $E_{\psi'}=E_{\phi'}\triangle E^0$. As a consequence, $E_\psi\sqcup E_{\psi'}=E_\phi\sqcup E_{\phi'}$ and $E_\psi\triangle E_{\psi'}=E_\phi\triangle E_{\phi'}$.
  \end{lemma}

\noindent(Indeed, by Lemma~\ref{lm:Ephi-Ephip}, each path $P\in \Pscr(\phi,\phi')$ is a concatenation of paths $Q_1,\ldots,Q_r$ such that consecutive $Q_j,Q_{j+1}$ are contained in different flows among $\phi,\phi'$ and either both leave or both enter their common vertex. Let $C_P:=C\triangle\,(\pi\cap Y)$ and $C'_P:=C'\triangle\,(\pi\cap Y')$. Exchanging the pieces $Q_j$ between $\phi$ and $\phi'$, we obtain an
$(XC_P|X'C'_P)$-flow $\alpha$ and an $(X\bar C_P|X'\bar{C'}_P)$-flow $\alpha'$ such that $E_{\alpha}=E_\phi\triangle\,E_{P}$ and $E_{\alpha'}=E_{\phi'}\triangle\,E_{P}$.

Doing so for all $P\in \Pscr(\phi,\phi')$ with $\pi(P)\in M_0$, we obtain flows
$\psi,\psi'$ satisfying the desired properties, taking into account that the
paths in $\Pscr(\phi,\phi')$ are pairwise disjoint. The uniqueness of
$\psi,\psi'$ is easy.)
\medskip

Lemma~\ref{lm:path_switch} implies that $M(\psi,\psi')=M(\phi,\phi')$ and $\Pscr(\psi,\psi')=\Pscr(\phi,\phi')$; therefore, the transformation of $\psi,\psi'$ by use of the paths in $\Pscr(\psi,\psi')$ related to $M_0$ returns $\phi,\phi'$.

 \begin{figure}[htb]
\begin{center}
\vspace{-0.1cm}
\includegraphics[scale=0.9]{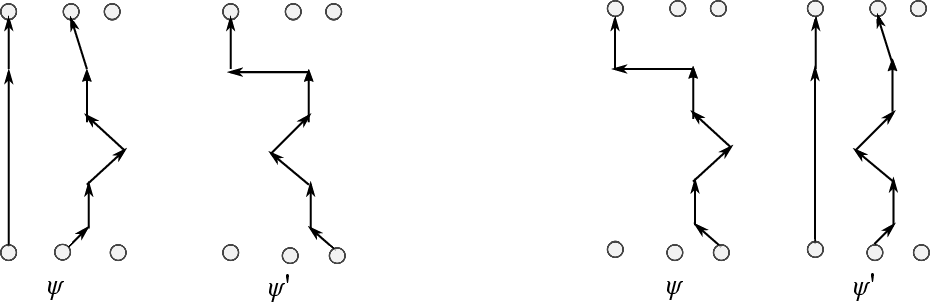}
\end{center}
\vspace{-0.5cm} 
\caption{Creating $(\psi,\psi')$ from $(\phi,\phi')$ in Fig.~\ref{fig:double}: \emph{left}: by use of $P$; \emph{right}: by use of $P'$}
 \label{fig:rearrange}
\end{figure}

Figure~\ref{fig:rearrange} illustrates flows $\psi,\psi'$ created from
$\phi,\phi'$ in Fig.~\ref{fig:double}. Here the left fragment shows
$\psi,\psi'$ when the exchange is performed by use of the single path
$P\in\Pscr(\phi,\phi')$ connecting $\hat s_2$ and $\hat s_3$,
and the right fragment shows $\psi,\psi'$ for the path $P'$ connecting 
$\hat s_1$ and $\hat t_2$ (see~(e) in Fig.~\ref{fig:double}).


\section{Proof of Theorem~\ref{tm:univ_relat}}  \label{sec:proof}

We first show implication (ii)$\to$(i) in Theorem~\ref{tm:univ_relat}. So let (ii) take place. For an $n\times n'$ TNN matrix $Q$, take a planar network $G=(V,E)$ and a weighting $w$ on $V$ such that the flow-generated function $f=f_{G,w}$ on $\Escr^{n,n'}$ defined by~\refeq{alg_f} satisfies the relation as in~\refeq{lind} for the given $Q$, i.e $\Delta_Q(I|J)=f(I|J)$ holds for all $(I,J)\in\Escr^{n,n'}$. For $X,Y,X',Y$ as above, consider families $\Ascr,\Bscr\Subset\CXYXY$ for which the matching relations as in (ii) of the theorem is valid. Then we argue as follows.

The summand in the R.H.S. of~\refeq{q_relat} concerning $(B,B')\in\Bscr$ can be expressed via double flows as follows:
  \begin{multline} \label{eq:ff-zeta}
  \Delta_Q(XB|X'B')\, \Delta_Q(X\bar B|X'\bar{B'})=f(XB|X'B')\, f(X\bar B|X'\bar{B'}) \\
=\left(\sum\nolimits_{\phi\in\Phi_{XB|X'B'}} w(\phi) \right) \times
 \left(\sum\nolimits_{\phi'\in\Phi_{X\bar B|X'\bar{B'}}} w(\phi') \right)  
                  \qquad\qquad\qquad\qquad \\
  =\sum\nolimits_{(\phi,\phi')\in\Dscr(B,B')} w(\phi)\, w(\phi')
  \qquad\qquad  \qquad\qquad \qquad\\
  =\sum\nolimits_{M\in\Mscr_{B,B'}}   \sum\nolimits_{(\phi,\phi')\in\Dscr(B,B')\,\colon M(\phi,\phi')=M} w(\phi)\,w(\phi'). \qquad
  \end{multline}
where for a flow $\phi''$, we write $w(\phi'')$ for $\prod(w(v)\colon v\in V_{\phi''}$), and denote by $\Dscr(B,B')$ the set of double flows $(\phi,\phi')$ for $(B,B')$ (i.e. $\phi$ is an $(XB|X'B')$-flow and $\phi'$ is an $(X\bar B|X\bar{B'})$-flow in $G$). The summand concerning $(A,A')\in\Ascr$ in the L.H.S. of~\refeq{q_relat} is expressed similarly.

Consider a matching $M\in\Mscr(\Bscr)$. Let $\Pi_M(\Bscr)$ be the family of proper pairs $(B,B')$ in $\Bscr$ for which $M$ is a feasible matching, i.e. $M\in \Mscr_{B,B'}$. Let $\Pi_M(\Ascr)$ be a similar family for $\Ascr$. Then $|\Pi_M(\Bscr)|=\#_M(\Bscr)$ and $|\Pi_M(\Ascr)|=\#_M(\Ascr)$. Therefore, the condition $\#_M(\Bscr)\le \#_M(\Ascr)$ (by (ii) in the theorem) implies that we can arrange an injective map $\sigma_M:\Pi_M(\Bscr)\to \Pi_M(\Ascr)$. 

Now for $(B,B')\in\Bscr$, consider a double flow $(\phi,\phi)\in\Dscr(B,B')$ and the corresponding matching $M=M(\phi,\phi')$; then $M\in\Mscr_{B,B'}$. Take the proper pair $(A,A')=\sigma_M (B,B')$. Then $M\in\Mscr_{A,A'}$, and therefore (cf.~\refeq{feas_match}(i)), $(A,A')$ is obtained from $(B,B')$ by the exchange operations w.r.t. some subset of couples $M_0\subseteq M$. Moreover, by transforming $(\phi,\phi')$ w.r.t. $M_0$ according to Lemma~\ref{lm:path_switch}, we obtain a double flow $(\psi,\psi')$ for $(A,A')$; denote it as $\tau(\phi,\phi')$. Since $E_\psi\sqcup E_{\psi'}=E_\phi\sqcup E_{\phi'}$ (again by Lemma~\ref{lm:path_switch}), we have $w(\psi)\, w(\psi')=w(\phi)\, w(\phi')$; so the contribution from the double flow $(\phi,\phi')$ to the $\Bscr$-sum of~\refeq{q_relat} is equal to that from $(\psi,\psi')$ to the $\Ascr$-sum. Also the map $\tau$ on the double flows $(\phi,\phi')$ concerning $\Bscr$ and related to $M$ is injective (taking into account that the transformation of $(\psi,\psi')$ as above w.r.t. $M_0$ returns $(\phi,\phi')$). 

It follows that the total contribution from the double flows concerning $\Bscr$ is equal to that from their images by $\tau$ concerning $\Ascr$. Now the desired inequality~\refeq{q_relat} is obtained by appealing to the expressions as in the last line of~\refeq{ff-zeta}, over all $(B,B')\in \Bscr$, and taking into account that all extra double flows for $\Ascr$ (that are not images by $\tau$) have nonnegative weights, in view of $w\ge 0$. 

(We also should take into account that, in the situation with $X=X'=\emptyset$, if some proper $(C,C')\in\CXYXY$ satisfies $C=C'=\emptyset$, then $\Delta_Q(XC|X'C')=1$ (by a convention in the beginning of Sect.~\SEC{theorem}), implying that $\Delta_Q(XC|X'C')\,\Delta_Q(X\bar C|X'\bar C')=\Delta_Q(Y|Y')$, which matches the agreement on flows that $f_{G,w}(\emptyset|\emptyset)=1$, whence $f_{G,w}(XC|X'C')\, f_{G,w}(X\bar C|X'\bar{C'})=f_{G,w}(Y|Y')$.) 

This completes the proof of (ii)$\to$(i) in the theorem.
  \smallskip

To obtain the converse implication (i)$\to$(ii), we use another result from~\cite{DKK2}. Given $X,Y,X',Y'$ as above, let $\Mscr_{X,Y,X',Y'}$ denote the set of feasible matchings occurring in $\Mscr(C,C')$ among all proper pairs $(C,C')\in\CXYXY$. The following assertion was shown in~\cite[pp.~461--66]{DKK2}.
 \begin{prop} \label{pr:P1P2}
For each $M\in\Mscr_{X,Y,X',Y'}$,  there exists a planar network $G=(V,E)$ with $n$ sources and $n'$ sinks possessing the following properties: for $(C,C')\in\CXYXY$,
 \begin{itemize}
\item[\rm(P1)] if $M\in\Mscr_{C,C'}$, then $G$ has a unique $(XC|\,X'C')$-flow and a unique $(X\bar C|\,X'\bar{C'})$-flow, i.e. $|\Phi^G_{XC|\,X'C'}|=|\Phi^G_{X\bar C|\,X'\bar{C'}}|=1$;
\item[\rm(P2)] if $M\notin\Mscr_{C,C'}$, then at least one set among $\Phi^G_{XC|\,X'C'}$ and $\Phi^G_{X\bar C|\,X'\bar{C'}}$ is empty.
  \end{itemize}
  \end{prop}
 
Relying on this, suppose that the family $\Bscr$ contains a matching $M$ violating condition~(ii) in the theorem, i.e. such that $\#_M(\Bscr)>\#_M(\Ascr)$. For this $M$, take a network $G=(V,E)$ as in Proposition~\ref{pr:P1P2} and assign unit weight $w(v)=1$ to each vertex $v\in V$. Then, by~(P1) and (P2), for the function $f=f_{G,w}$ and a proper pair $(C,C')\in\CXYXY$, both values $f(XC|\,X'C')$ and $f(X\bar C|\,X'\bar{C'})$ are ones if $M\in\Mscr(C,C')$, and at least one of them is zero otherwise. This implies that the $n\times n'$ matrix $Q$ corresponding to $f_{G,w}$  (see~\refeq{lind}) is totally nonnegative and satisfies
  \begin{multline*}
  \sum\nolimits_{(A,A')\in \Ascr} \Delta_Q(XA|\,X'A')\, \Delta_Q(X\bar A|\,X'\bar A') 
  =\#_M(\Ascr)\\
      <\#_M(\Bscr)=\sum\nolimits_{(B,B')\in \Bscr} \Delta_Q(XB|\,X'B')\, \Delta_Q(X\bar B|\,X'\bar B'). \qquad
   \end{multline*}

So~\refeq{q_relat} is violated. This yields (i)$\to$(ii), completing the proof of Theorem~\ref{tm:univ_relat}. 
\smallskip

Note that the assertion as in the above proposition is proved in~\cite{DKK2} by use of an explicit construction that, given a matching $M\in \Mscr_{X,Y,X',Y'}$, devises a graph $G$ as required; the size of this graph is smaller than $(n+n')^2$. 

\begin{figure}[htb]
\begin{center}
\vspace{-0.1cm}
\includegraphics[scale=0.8]{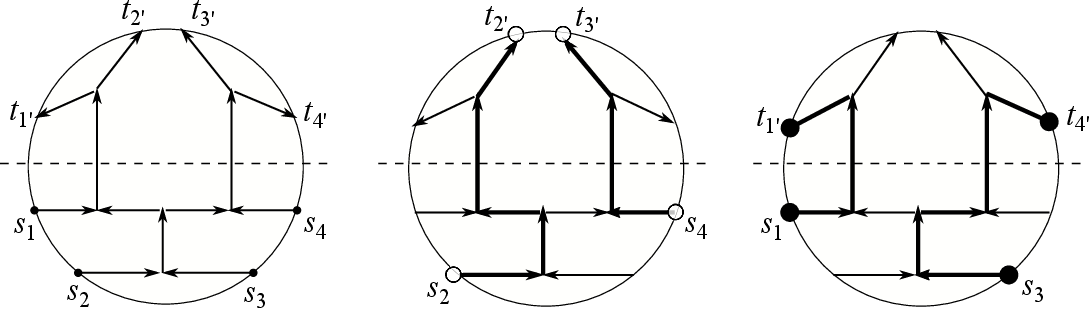}
\end{center}
\vspace{-0.5cm} 
\caption{\emph{left}: graph $G$; \;\;\emph{middle}: $(24|2'3')$-flow; \;\;\emph{right}: $(13|1'4')$-flow}
 \label{fig:graphG}
 \end{figure}

A simple example of $G$ is drawn in the left fragment of Fig.~\ref{fig:graphG}. Here $X=X'=\emptyset$,  $Y=[4]$, $Y'=\{1',\ldots,4'\}$, and $M$ consists of the couples $14,\, 23,\, 1'2',\, 3'4'$. The middle and right fragments illustrate, respectively, the unique $(C|C')$-flow and unique $(\bar C|\bar{C'})$-flow (drawn in bold), where $(C,C')=(\{2,4\},\{2',3'\})$ (satisfying $M\in\Mscr_{C,C'}$).

To finish this section, note that, due to the criterion in Theorem~\ref{tm:univ_relat}, the task of recognizing whether a given pair $(\Ascr,\Bscr)$ provides a universal inequality or not looks rather straightforward, as it reduces to computing the set $\Mscr_{C,C'}$ of feasible matchings for each $(C,C)\in\Ascr\cup\Bscr$ (which is  routine, though $|\Mscr_{C,C'}|$ may be exponentially large in $|Y|+|Y'|$), followed by comparing the numbers $\#_M(\Ascr)$ and $\#_M(\Bscr)$. 


\section{Illustrations}  \label{sec:illustr}

Here we give two appealing examples to illustrate how the matching approach works.

\subsection{Example 1} \label{ssec:ex1}
Let $Q$ be an $n\times n$ TNN matrix. As before, we will denote column indexes with primes. Consider subsets $A,B\subset[n]$ of the same size: $|A|=|B|$. We assert that the relation (of submodular character)
  \begin{equation} \label{eq:adam}
  \Delta_Q(A|A')\,\Delta_Q(B|B')\ge \Delta_Q((A\cup B)|(A'\cup B'))\, 
                 \Delta_Q((A\cap B)|(A'\cap B'))
  \end{equation}
 is valid, i.e. \refeq{adam} gives an instance of universal inequalities for TNN matrices. Here for $C\subseteq[n]$, we write $C'$ for $\{i'\colon i\in C\}$. 

Denote $A\cap B$ by $X$; $A-B$ by $\tilde A$; $B-A$ by $\tilde B$; and $\tilde A\cup\tilde B$ by $Y$. Then~\refeq{adam} can be rewritten as
  \begin{equation} \label{eq:adamp}
\Delta_Q(\tilde A X|\tilde A' X')\,\Delta_Q(\tilde B X|\tilde B'X')
    \ge \Delta_Q(XY|X'Y')\, \Delta_Q(X|X').
    \end{equation}

Then the family $\Ascr$ and $\Bscr$ consist of single pairs in $\CXYXY$, namely, $(\tilde A,\tilde A)$ and $(Y,Y')$, respectively. To show~\refeq{adamp}, refer to the elements in $\tilde A,\tilde A', Y, Y'$ as white, and in $\tilde B,\tilde B$ as black. Then the (proper) pair $(Y,Y')$ occurring in the R.H.S. of~\refeq{adamp} admits only one feasible matching $M$; namely, $M$ consists of the couples $ii'$ for all $i\in Y$ (which connect the white elements in the lower half of circumference $O$ with their white copies in the upper half). Since each of these couples $ii'$ connects elements of the same color in the diagram for the pair $(\tilde A,\tilde A')$ in the L.H.S. of~\refeq{adamp} (namely, connecting white elements $i,i'$ when $i\in \tilde A$, and black elements when $i\in\tilde B$), the matching $M$ on $Y\sqcup Y'$ is feasible for the pair $(\tilde A,\tilde A')$ as well. This proves~\refeq{adamp} (and~\refeq{adam}) by Theorem~\ref{tm:univ_relat}.

(By the way, in the case $A\cap B=\emptyset$, the R.H.S. in~\refeq{adam} reduces to $\Delta_Q((A\cup B)|(A'\cup B'))$, taking into account that $\Delta_Q(\emptyset|\emptyset)=1$.)

In addition, note that $(\tilde A,\tilde A')$ admits more feasible matchings; one of them is
  \begin{numitem1} \label{eq:Mp}
 $M'$ formed by the ``horizontal'' couples $\{\alpha_i,\beta_{k-i+1}\}$ and $\{\alpha'_i,\beta'_{k-i+1}\}$, $i=1,\ldots,k=|\tilde A|$, where $\tilde A=(\alpha_1<\cdots<\alpha_k)$ and $\tilde B=\{\beta_1<\cdots<\beta_k)$.
  \end{numitem1}
  
Note also that, instead of principal minors, one can take $A'=\sigma(A)$ and $B'=\sigma(B)$, where $\sigma$ is an order preserving injective map of $A\cup B$ into $[n]$.

\subsection{Example 2} \label{ssec:ex2}
Let $A,B,\tilde A,\tilde B,X,Y$ be as in the previous example. We can try to sharpen expression~\refeq{adam} as follows:
  \begin{multline} \label{eq:dodg}
  \Delta_Q(A|A')\,\Delta_Q(B|B') \\\ge \Delta_Q((A\cup B)|(A'\cup B'))\, 
                 \Delta_Q((A\cap B)|(A'\cap B')) + \Delta_Q(A|B')\,\Delta_Q(B|A').
  \end{multline}
This can be rewritten as (cf.~\refeq{adamp}
   \begin{multline} \label{eq:dodgp}
\Delta_Q(\tilde A X|\tilde A' X')\,\Delta_Q(\tilde B X|\tilde B'X') \\
    \ge \Delta_Q(XY|X'Y')\, \Delta_Q(X|X') + \Delta_Q(\tilde A X|\tilde B' X')\,   \Delta_Q(\tilde BX|\tilde A'X').
    \end{multline}
 
So the family $\Ascr$ consists of the single proper pair $(\tilde A,\tilde A)$ in $\CXYXY$, and $\Bscr$ of two proper pairs $(Y,Y')$ and $(\tilde A,\tilde B')$. In the particular case when $|\tilde A|=|\tilde B|=1$, this turns into equality for any real matrix $Q$, which is known as \emph{Dodgson condensation formula} (going back to~\cite{dodg}). For this reason, we conditionally refer to~\refeq{dodg} (and~\refeq{dodgp}) in a general case concerning TNN matrices $Q$ as a \emph{generalized inequality of Dodgson's type}. Note that it need not hold in general, and it is tempting to characterize the pairs $(A,B)$ (with $|A|=|B|$) for which~\refeq{dodg} is valid for all TNN matrices. 

Recall that sets $A,B\subseteq[n]$ are called \emph{strongly separated} if there are no triple $i<j<k$ such that $i,k$ belong to one, and $j$ to the other set among $A-B$ and $B-A$. When $|A|=|B|$, the sets are called \emph{weakly separated} if there are no quadruple $i<j<k<\ell$ such that $i,k$ belong to one, and $j,\ell$ to the other set among $A-B$ and $B-A$ (the definition is somewhat different when $|A|\ne|B|$, but this is not needed to us). We show the following
 \begin{prop} \label{pr:gen_dodg}
Let $A,B\subset[n]$ be weakly separated and $|A|=|B|$. The following are equivalent:
  
  {\rm(i)} \refeq{dodg} is valid for every $n\times n$ TNN matrix $Q$;
    
  {\rm(ii)} $A$ and $B$ are strongly separated.
    \end{prop}
    \begin{proof}
(ii)$\to$(i).
For $A,B$ as in (ii), one may assume that $\tilde A=(\alpha_1<\cdots<\alpha_k)$, $\tilde B=(\beta_1<\cdots<\beta_k)$ and $\alpha_k<\beta_1$. Let $M$ and $M'$ be as defined in Sect.~\SSEC{ex1}. Then $M$ is the unique feasible matching for the proper pair $(Y,Y')$ and is a feasible matching for the proper pair $(\tilde A,\tilde A')$. A similar property takes place for the other pair $(\tilde A,\tilde B')$ occurring in the R.H.S. of~\refeq{dodgp}. More precisely, one can see that this pair admits a unique feasible matching; this is just $M'$ defined in~\refeq{Mp}. At the same time, $M'$ is feasible for $(\tilde A,\tilde A')$ as well.
So the desired property~(i) follows from Theorem~\ref{tm:univ_relat}.

The picture below illustrated the matching $M'$ for $(\tilde A,\tilde A')$ (left) and for $(\tilde A,\tilde B')$ (right); here $\tilde A=\{1,2\}$ and $\tilde B=\{3,4\}$.

\vspace{-0cm}
\begin{center}
\includegraphics[scale=0.75]{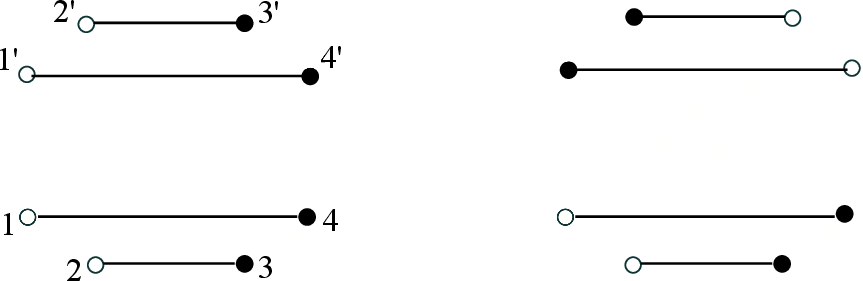}
\end{center}
\vspace{-0cm}

(i)$\to$(ii).
Now let $A,B\subset[n]$ with $|A|=|B|$ be weakly but not strongly separated. W.l.o.g., we may assume that $X=A\cap B=\emptyset$ (since feasible matchings do not depend on $X$) and that $Y=\tilde A\cup\tilde B=A\cup B=[2k]$. Also assume that $A$ contains 1.
  \smallskip

First suppose that $k=2$. Then there is only one case, namely, $A=\{1,4\}$ and $B=\{2,3\}$. This case is illustrated in the picture.

\vspace{-0cm}
\begin{center}
\includegraphics[scale=0.75]{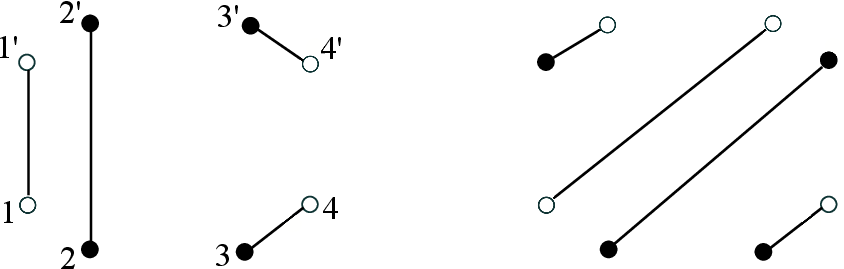}
\end{center}
\vspace{-0cm}

Here the left (right) fragment concerns the proper pair $(A,A')$ (resp. $(A,B')$). One can see that the matching drawn in the right fragment is feasible for $(A,B')$ but not for $(A,A')$. In turn, the matching drawn in the left fragment is feasible for $(A,A')$ but not for $(A,B')$. Therefore, for the collections $\Ascr=\{(A,A')\}$ and $\Bscr=\{(Y,Y'),(A,B')\}$, the sets $\Mscr(\Ascr)$ and $\Mscr(\Bscr)$ are incomparable (i.e. neither $\Mscr(\Ascr)\supseteq\Mscr(\Bscr)$ nor $\Mscr(\Ascr)\subseteq\Mscr(\Bscr)$ takes place). Then, by Theorem~\ref{tm:univ_relat}, neither~\refeq{dodgp} nor the  reversed inequality to it is universal for TNN matrices.

Next suppose that $k\ge 3$. Since $A,B$ are weakly but not strongly separated and $1\in A$, there are $p,q\in[2k]$ such that $A=\{1,\ldots,p\}\cup \{q+1,\ldots,2k\}$ and $B=\{p+1,\ldots, q\}$. Also $q-p=|B|=k\ge 3$. This implies at least one of the following: $p\ge 2$, or $2k-q\ge 2$. One may assume the former (the latter is symmetric). It is easy to see that: (a) the pair $\{p,p+1\}$ has one element in $A$ and the other in $B$, and similarly, the pair $\{p',(p+1)'\}$ has one element in $A'$ and the other in $B'$; and (b) the sets $\hat A:=A-\{p\}$ and $\hat B:=B-\{p+1\}$ are again weakly but not strongly separated. Also $|\hat A|=|\hat B|=k-1\ge 2$. The by induction we may assume the existence of matchings $\hat M, \hat M'$ on $(Y-\{p,p+1\})\sqcup (Y-\{p',(p+1)'\})$ such that $\hat M$ is feasible for $(\hat A,\hat A')$ but not for $(\hat A,\hat B')$, whereas the behavior of $\hat M'$ is converse. Now adding to each of $\hat M,\hat M'$ the ``horizontal'' couples $\{p,p+1\}$ and $\{p',(p+1)'\}$, we obtain the desired matching for $(A,A')$ and $(A,B')$, showing that neither~\refeq{dodgp} nor its reversed inequality is universal for TNN matrices for $A,B$ in question.

This completes the proof of the proposition.   
    \end{proof}
    
\noindent\textbf{Remark.} When sets $A,B$ are not weakly separated, the question whether~\refeq{dodg} is universal or not becomes more intricate. On the one hand, one can easily extend the construction from the proof of (i)$\to$(ii) in Proposition~\ref{pr:gen_dodg} to obtain a plenty of instances when~\refeq{dodg} need not hold. On the other hand, there are cases for which it is always valid. This is so already for the simplest example $A=\{1,3\}$ and $B=\{2,4\}$ illustrated in the following picture. 

\vspace{-0cm}
\begin{center}
\includegraphics[scale=0.75]{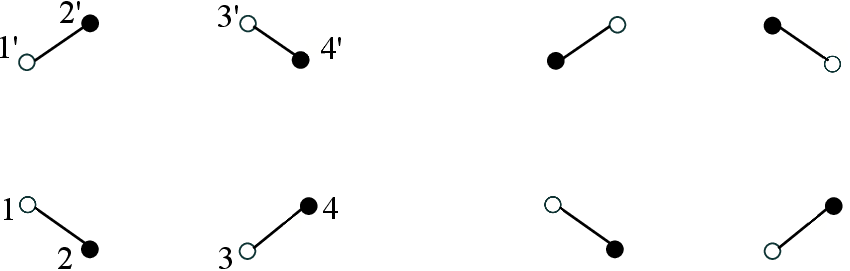}
\end{center}
\vspace{-0cm}

Here the drawn matching is feasible for both $(A,B)$ and $(A,B')$. One can see that $(A,B')$ has no feasible matchings containing a ``vertical'' couple (connecting $Y$ and $Y'$) and that any feasible matching formed by ``horizontal'' couples is feasible for $(A,B)$ as well (there are 4 such matchings in all). Therefore,~\refeq{dodg} is valid, yielding a TNN-universal instance.

\end{document}